\documentclass{ifacconf}

\usepackage{natbib}            
\usepackage{graphicx}          

\usepackage[utf8]{inputenc}
\usepackage[T1]{fontenc}

\usepackage[retainorgcmds]{IEEEtrantools}

\usepackage{amssymb,amsmath}
\usepackage{url}

\graphicspath{{Figures/}}
\usepackage{subfigure}

\begin{document}

\begin{frontmatter}

\title{Controllability of the 1D Schr\"odinger equation by the flatness approach}

\author[First]{Philippe Martin}
\author[Second]{Lionel Rosier}
\author[First]{Pierre Rouchon}

\address[First]{Centre Automatique et Syst\`emes, MINES ParisTech, 75272 Paris, France (e-mail: $\{$philippe.martin,pierre.rouchon$\}$ @mines-paristech.fr).}
\address[Second]{Institut Elie Cartan, UMR 7502 UdL/CNRS/INRIA, BP 70239, 54506 Vand\oe uvre-l\`es-Nancy, France (e-mail: Lionel.Rosier@univ-lorraine.fr)}

\begin{keyword}                           
Partial differential equations, Schr\"odinger equation, beam equation, boundary control, exact controllability, path planning, flatness.
\end{keyword}                             

\begin{abstract}                          
We derive in a straightforward way the exact controllability of the 1-D Schr\"odinger equation with a Dirichlet boundary control. We use the so-called {\em flatness approach},
which consists in parameterizing the solution and the control by the derivatives of a ``flat output''. This provides an explicit control input achieving the exact controllability in the
energy space. As an application, we derive an explicit  pair of  control inputs achieving the exact steering to zero for a simply-supported beam.
\end{abstract}

\end{frontmatter}

\section{Introduction} The exact controllability of the linear Schr\"odinger equation (or of the plate equation) was investigated in
\cite{lions-book}, \cite{machtyngier}, \cite{komornik-book} with the multiplier method,
in  \cite{haraux}, \cite{jaffard}, \cite{KL-book} with nonharmonic Fourier analysis, in \cite{lebeau} with microlocal analysis, and in \cite{CFNS,liu} with frequency domain tests.
The exact controllability was extended to the semilinear Schr\"odinger
equation in \cite{RZ2009S1D,RZ2009SND,RZ2010},
and \cite{laurentS1D,laurentS3D} by means of Strichartz estimates and Bourgain analysis.

All the above results rely on some observability inequalities for the adjoint system. A direct approach which does not involve the adjoint
system was proposed in \cite{LM,R2002,LT}. The result in \cite{R2002} used a fundamental solution of the Schr\"odinger equation
with compact support in time, and provided controls that are Gevrey.

In this paper, we derive in a straightforward way the null (or, equivalently, exact) controllability of the (linear) Schr\"odinger
equation on the interval $(0,1)$ with a Dirichlet control at $x=1$. More precisely, for any final time $T>0$ and
any $\theta _0\in L^2(0,1)$, we provide an explicit and regular control such that the state reached at time $T$ is exactly zero.
We use the so-called {\em flatness} approach \cite{FLMR}, which consists in parameterizing the solution $\theta$ and the control $u$
by the derivatives of a ``flat output'' $y$; this notion was initially introduced for finite-dimensional (nonlinear) systems,
and later extended to PDEs, \cite{LMR, LR,MTK,meurer-book}. In \cite{MRR1,MRR2,MRR3}, we proved that the null controllability
of the heat equation could be derived with the flatness approach, and that this approach provided very efficient numerical schemes by taking
partial sums in the series.
The design of the control in \cite{MRR1} was done in two steps. In the first one, a
null control was used to reach an intermediate state Gevrey of order $1/2$, while a flatness-based control was used in the second step to drive this (regular)
intermediate state to 0. Clearly, this strategy cannot be used for the Schr\"odinger equation, as an application of a null control in the first step would yield a
state which may be only $L^2$. Instead, using some idea in \cite{LT,RZ2009SND}, we apply here a {\em nonull} control in the first step, using a smoothing
effect for the Schr\"odinger equation on $\mathbb R$, to reach an intermediate state which is again Gevrey of order $1/2$. The second step is then carried out as in \cite{MRR1}.

The paper is outlined as follows. In Section 2 we consider the control problem for the Schr\"odinger equation in dimension $N=1$.  In Proposition
\ref{prop1} we investigate an ill-posed problem with Cauchy data in a Gevrey class and prove its (global) well-posedness. Theorem \ref{thm1} then establishes
the null controllability in small time for any initial data in $L^2$. Section 3 is devoted to the exact controllability of a simply-supported beam.
 In Theorem \ref{thm2} we derive a null controllability result by the flatness approach for the beam equation with two boundary controls. Some  additional computations
 needed to compute numerically the control for the Schr\"odinger equation are provided in Section 4.

Let us introduce some definitions and notations. We say a function $y\in C^\infty([0,T])$ is {\em Gevrey of order $s\ge 0$ on $[0,T]$}
if there exist some positive constants $M,R$ such that
\begin{equation}
\label{G1}
|y^{(p)}(t)| \le M\frac{ (p!)^s}{R^p} \qquad \forall t\in [0,T],\ \forall p\ge 0;
\end{equation}
If $y$ is complex-valued the definition applies for the real and complex parts.
More generally, for any compact set $K\subset {\mathbb R}^N$ ($N\ge 1$), and any function $y:x=(x_1,x_2,...,x_N)\in K\mapsto y(x)\in {\mathbb R}$ which is
of class $C^\infty$ on $K$ (i.e. $y$ is  the restriction to $K$ of a function of class $C^\infty$ on some open neighborhood $\Omega$ of $K$),
we shall say that {\em $y$ is Gevrey of order $s_1$ in $x_1$, $s_2$ in $x_2$, ..., $s_N$ in $x_N$ on $K$},
where $s_i\ge 0$ for $1\le i\le N$, if there
are some positive constants $M,R_1,...,R_N$ such that
\begin{equation}
\label{G3}
|\partial _{x_1}^{p_1}\partial_{x_2}^{p_2}\cdots \partial_{x_N}^{p_N}y(x)| \le
M\frac{\prod_{i=1}^N (p_i !)^{s_i}}{\prod_{i=1}^N R_i^{p_i} }, \quad \forall x\in K,\ \forall p\in {\mathbb N} ^N.
\end{equation}
\section{Controllability of the Schr\"odinger equation}
We are concerned with the (null) controllability of the system
\begin{IEEEeqnarray}{rCll}
i\theta _t +\theta _{xx} &=&0, \quad &(t,x)\in (0,T)\times (0,1), \label{B1}\\
\theta (t,0) &=& 0, \quad &t\in (0,T), \label{B2}\\
\theta (t,1) &=& u(t),\quad &t\in (0,T), \label{B3}\\
\theta (0,x) &=& \theta _0 (x), \quad &x\in (0,1), \label{B4}
\end{IEEEeqnarray}
where $\theta _0\in L^2(0,1)$ is given; $\theta,\theta_0$ and~$u$ are complex-valued functions.

Let
\begin{equation}
\label{AA5}
y(t)=\theta _x(t,0),\qquad t\in [0,T].
\end{equation}
We claim that \eqref{B1}-\eqref{B2} is ``flat'' with $y$ as flat output, which means that the map $\theta \to y$ is a bijection
between  appropriate spaces of smooth functions. 
Indeed, let us seek a formal solution $(\theta,u)$  of \eqref{B1}-\eqref{B3} and \eqref{AA5} in the form
\[
\theta (t,x) = \sum_{k\ge 0} \frac{x^k}{k!} a_k(t), \quad u(t)= \sum_{k\ge 0} \frac{a_k(t) }{k ! },
\]
 where the $a_k$'s are functions yet to define. Plugging the formal solution into \eqref{B1}  yields
 \[
 \sum_{i\ge 0} \frac{x^k}{k!} [a_{k+2} +i a_k'] =0.
 \]
 Thus
 \[
 a_{k+2}= -i a_k', \quad \forall k\ge 0.
 \]
 On the other hand, we infer  from \eqref{B2} and  \eqref{AA5}  that $a_0(t)=0$ and $a_1(t)=y(t)$. It follows that for all $k\ge 0$
 \[
 a_{2k}=0,\quad a_{2k+1} = (-i)^k y^{(k)}.
 \]
Eventually,
\begin{IEEEeqnarray}{rCl}
\theta (t,x) &=& \sum_{k\ge 0} \frac{x^{2k +1 } }{ (2k +1) ! } (-i)^k y^{(k)} (t) \label{AA10a}\\
 u(t) &=&\sum _{k\ge 0} \frac{ (-i)^k y^{(k)} (t)}{(2k+1) !} \cdot \label{AA10b}
\end{IEEEeqnarray}
In particular $\theta$ is uniquely defined in terms of $y$.

 Conversely, if $y\in C^\infty ( [0,T] )$ and the first formal series in \eqref{AA10a} is convergent in $C^2([0,T]\times [0,1])$, then it is easily seen that
 \eqref{B1}-\eqref{B3} and \eqref{AA5} hold.
 Our first result shows that the above computations are fully justified when $y$ is Gevrey of order $s\in [0,2)$.

 \begin{prop}
 \label{prop1}
 Let $s\in [0,2)$, $-\infty <t_1<t_2<\infty$, and $y\in C^\infty ([t_1,t_2] )$ satisfying for some constants $M,R>0$
 \begin{equation}
 \label{AA11}
 |y^{(k)}(t) | \le M \frac{k!^s}{R^k}, \qquad \forall k\ge 0,\ \forall t\in [t_1,t_2].
 \end{equation}
 Then the function $\theta$ defined in \eqref{AA10a} is Gevrey of order $s$ in $t$ and $s/2$ in $x$ on
 $[t_1,t_2]\times [0,1]$.
 \end{prop}
 \begin{pf}
 We want to prove the formal series
 \begin{equation}
 \label{AA12}
 \partial _t ^m\partial _x ^n \theta (t,x) = \sum_{2k+1 \ge n}  \frac{ x^{2k+1-n } }{ (2k + 1 -n)! } (-i)^k y^{(k+m)}(t)
 \end{equation}
 is uniformly convergent on $[t_1,t_2]\times [0,1]$  with an estimate of its sum of the form
 \begin{equation}
 \label{AA12bis}
|\partial _t ^m\partial _x ^n \theta (t,x)| \le C
\frac{m! ^s}{R_1^m}\, \frac{n! ^\frac{s}{2}}{R_2^n} \cdot
 \end{equation}
By \eqref{AA11}, we have for all $(t,x)\in [t_1,t_2]\times [0,1]$
\begin{IEEEeqnarray*}{rCl}
&&\left\vert  \frac{x^{2k+1 -n}}{(2k+1-n)!} y^{(k+m)} (t) \right\vert \\
&&\quad   \le\frac{M}{R^{k+m}} \, \frac{(k+m)! ^s}{(2k + 1-n)!} \\
&&\quad \le \frac{M}{R^{k+m} }\,  \frac{(2^{k+m} k! \, m! )^s}{(2k +1 -n)!} \\
&&\quad \le \frac{M}{R_1^{k+m}}\, \frac{(2^{-2k} \sqrt{\pi k} \, (2k)! )^\frac{s}{2} }{(2k + 1-n)!}\, m!^s \\
&&\quad \le \frac{M}{R_1^{k+m}}\, \frac{(2^{-2k-1} \sqrt{\pi k} \, (2k+1)! )^\frac{s}{2} }{(k+\frac{1}{2})^\frac{s}{2} (2k + 1-n)!}\, m!^s \\
&&\quad  \le  M \frac{(\pi k) ^{\frac{s}{4}} }{R_1 ^k (k+\frac{1}{2})^\frac{s}{2}  (2k +1 -n)! ^{1-\frac{s}{2}} } n! ^{\frac{s}{2}} \frac{m!^s}{R_1 ^m},
\end{IEEEeqnarray*}
where we have set $R_1 = 2^{-s}R$ and used twice\\ $(p+q)! \le 2^{p+q} p! q! $, and
\begin{equation}
\label{AA13}
(2k) ! \sim \frac{2^{2k} }{\sqrt{\pi k} } k! ^2
\end{equation}
which follows at once from Stirling's formula. Since $\sum_{2k+1\ge n} \frac{(\pi  k)^\frac{s}{4} }{R_1 ^k  (k+\frac{1}{2})^\frac{s}{2}  (2k+1-n)! ^{1-\frac{s}{2} }  } <\infty$,
we infer the uniform convergence of the series in \eqref{AA12} for all $m,n\ge 0$. This shows that $\theta \in C^\infty ([t_1,t_2]\times [0,1])$.
On the other hand, picking any $R_2\in (0,\sqrt{R_1})$, since
\begin{IEEEeqnarray*}{rCl}
&&M\sum_{2k +1\ge n} \frac{ (\pi  k)^\frac{s}{4} }{R_1 ^k   (k+\frac{1}{2})^\frac{s}{2}  (2k+1-n)! ^{1-\frac{s}{2} }  } \\
&&\quad \le  K R_1 ^{-\frac{n}{2}}  \sum_{j\ge 0} \frac{ j^\frac{s}{4} + n^\frac{s}{4} }{R_1 ^\frac{j}{2}  j! ^{1-\frac{s}{2} }  }\\
&&\quad \le C R_2^{-n},
\end{IEEEeqnarray*}
for some constants $K,C>0$ independent of $n$,
we conclude that
\[
|\partial _t ^m \partial _x^n \theta (t,x)| \le
C \frac{n! ^{\frac{s}{2}} }{R_2^n} \, \frac{m! ^s}{R_1^m},
\]
which proves that $\theta$ is Gevrey of order $s$ in $t$ and $s/2$ in $x$, as desired.\qed

 \end{pf}

 Let $\theta _0\in L^2(0,1)$  be given. Define $v_0\in L^2({\mathbb R} )$ as
 \begin{equation}
 v_0(x) =
 \left\{
 \begin{array}{ll}
 \theta  _0(x)\quad  &\textrm{ if } x\in (0,1),\\
 -\theta _0(-x) &\textrm { if } x\in (-1,0), \\
 0 &\textrm{ if } x\in (-\infty , -1)\cup (1, +\infty ).
 \end{array}
 \right.
 \label{v0}
 \end{equation}
 Let $v=v(t,x)$ denote the solution of the Cauchy problem
 \begin{IEEEeqnarray}{rCl}
 iv_t + v_{xx} =0,&&\quad (t,x)\in {\mathbb R} ^2,  \label{v1}\\
 v(0,x)=v_0(x), && \quad x\in {\mathbb R}.  \label{v2}
 \end{IEEEeqnarray}
The following properties are well known, see e.g. \cite{cazenave-book,LP-book,RZ2009SND}:
\begin{IEEEeqnarray}{rCl}
&&v(t,-x) = - v(t,x)\quad \textrm{for a.e. } x \in {\mathbb R}, \textrm{ for all } t\in {\mathbb R},  \label{P0} \\
&&v\in C^\infty ( {\mathbb R} _t \setminus \{0\} \times {\mathbb R} _x ), \label{P1}\\
&&\left( \int_{-\infty} ^\infty || v ||^4_{L^\infty ({\mathbb R}) } dt \right)^{\frac{1}{4}} \le c ||v_0||_{L^2(\mathbb R  )}, \label{P2} \\
&&\sup_x \int_{-\infty}^\infty |D_x ^{\frac{1}{2}}v|^2dt \le c ||v_0||_{L^2(\mathbb R ) }^2. \label{P3}
\end{IEEEeqnarray}
\eqref{P1}  rests on the fact that $v_0$ is {\em compactly supported}. To justify \eqref{P1}, it is sufficient to introduce the operator
\[ Pu:=(x+2it\partial _x)u = 2it e^{i\frac{|x|^2}{4t}} \partial _x \big( e^{- i\frac{|x|^2}{4t}} u  \big)  \]
and to notice that it commutes with the Schr\"odinger operator $L= i\partial _t + \partial _{xx}$. The same is true for  $P^k$ and $L$ for all $k\in \mathbb N$, and hence
\[ || (P^k v)(t,.) ||_{ L^2 (\mathbb R ) } =   || (P^k v) (0,.)||_{ L^2(\mathbb R ) } = ||x^k v_0||_{L^2(\mathbb R )} <\infty , \]
where we used the conservation of the $L^2(\mathbb R)-$norm for the solutions of \eqref{v1}-\eqref{v2}.
This yields $v\in L^\infty ( (t_1,t_2), H^k_{loc}( \mathbb R ))$ for all $0<t_1<t_2$ and $k\in \mathbb N$, and hence \eqref{P1} by using  \eqref{v1} inductively.
 In particular, it follows from \eqref{P1} and \eqref{P2} that for all $\tau >0$,
 \begin{equation}
 \label{P4}
 v(.,1)\in C^\infty ((0,\tau])\cap L^4(0,\tau ).
 \end{equation}

 With Proposition \ref{prop1} at hand, we can derive a null controllability result obtained in a constructive way by the flatness approach.

\begin{thm}
\label{thm1}
Let $\theta _0\in L^2(0,1)$ and $T>0$ be given. Let $v$ denote the solution of \eqref{v1}-\eqref{v2}, where $v_0$ is as in \eqref{v0}.
Pick any $\tau \in (2T/3,T)$ and any $s\in (1,2)$. Then there exists a function $y:[\tau , T]\to \mathbb R$
Gevrey of order $s$ on $[\tau , T]$ such that, setting
 \begin{equation}\label{eq:u}
u(t) = \left\{
\begin{array}{ll}
v(t,1)\quad &\text{if }\  0< t\le \tau ,\\
\sum _{k\ge 1} \frac{ (-i)^k y^{(k)} (t)}{(2k+1) !}  \quad &\text{if }\  \tau  < t \le T,
\end{array}
\right.
\end{equation}
the solution $\theta$ of \eqref{B1}-\eqref{B4} satisfies $\theta (T,.) = 0$.  Furthermore, the control function $u$ is in $L^4(0,T)$ and it is Gevrey of order
$s$ in $t$ on $[ \varepsilon ,T]$ for all $\varepsilon \in (0,T)$,
$\theta \in C([0,T],L^2(0,1)) \cap C^\infty ((0,T]\times [0,1])$, and $\theta$ is Gevrey of order $s$ in $t$ and $s/2$ in $x$ on
$[\varepsilon ,T]\times [0,1]$ for all $\varepsilon \in (0,T)$.
\end{thm}
\begin{pf} The idea is to apply first a control in the time interval $[0,\tau ]$ to smooth out the state function, and next to use the above flatness
approach to steer the (more regular) state function to 0 in the time interval  $[\tau , T]$.\\
{\sc Step 1. Free evolution.}\\
We set $u(t)=v(t,1)$ for $t\in [0,\tau ]$. By~\eqref{P0}-\eqref{P1} $v$ is smooth and odd for all $t>0$ hence $v(t,0)=0$.
Therefore, $\theta (t,x)=v(t,x)$ for $(t,x)\in [0,\tau ] \times [0,1]$.
Introduce the usual fundamental solution of the Schr\"odinger equation
\[
E(t,x)=\frac{1}{ ( 4\pi i t )^\frac{1}{2} } e^{ i \frac{x^2}{4t} }.
\]
Then
\begin{equation}
v(t,x)=(E(t,.) * v_0)(x) =\int_{-1}^1 E(t,x-y)v_0(y)dy.
\label{C100}
\end{equation}
\begin{lem}
\label{lem1}
The function $v(t,x)$ is Gevrey of order 1 in $t$ and $1/2$ in $x$ on $[t_1,t_2]\times [-L,L]$ for all $0<t_1<t_2$ and $L>0$. Furthermore, we can write
\begin{equation}
v(\tau , x) = \sum_{k\ge 0} y_k (-i) ^k \frac{x^{2k+1 }}{ (2k +1 ) !}, \quad x\in \mathbb R \label{B110}
\end{equation}
with
\begin{equation}
|y_k | \le C (\tau )||\theta _0||_{ L^1(0,1) }\Bigl( \frac{2}{\tau } \Bigr) ^k  k !
\label{B12}
\end{equation}
where $C(\tau )$ is a continuous function on $(0,+\infty )$.
\end{lem}
\noindent
{\em Proof of Lemma \ref{lem1}:}
Pick any $k\in \mathbb N$. Then
\begin{IEEEeqnarray}{rCl}
\partial _x ^k v(t,x) &=& \int_{-1}^1 (\partial _x ^k E) (t,x-y) v_0(y)\, dy  \nonumber \\
&=& \frac{1}{ ( 4\pi i t )^\frac{1}{2} }\int_{-1}^1 \partial_x [ e^{ i \frac{x^2}{4t} } ] \label{W3}
(t,x-y) v_0(y) dy.
\end{IEEEeqnarray}
Let us first check that $e^{ i \frac{x^2}{4t} }$ is Gevrey of order $1/2$ in $x$. Clearly
\[
e^{ i \frac{x^2}{4t} }=\sum_{k\ge 0} \Bigl(\frac{i}{4t}\Bigr)^k \frac{x^{2k}}{k!} =\sum _{l\ge 0 } a_l \frac{x^l }{ l ! } ,
\]
where
\[
a_l = \left\{
\begin{array}{ll}
\Bigl(\frac{i}{4t}\Bigr)^k \frac{(2k)!}{k!} & \textrm{ if } l=2k,\\
0 &\textrm{ if } l=2k+1.
\end{array}
 \right.
\]
Then for $t_1\le t\le t_2$ and $x\in [-L,L]$,
\begin{IEEEeqnarray}{rCl}
|a_{2k} |  &\le & \frac{c}{(4t)^k} \frac{ (2k) !} {(2k) ! ^\frac{1}{2} 2^{-k}  (\pi k ) ^{\frac{1}{4}} } \nonumber \\
&\le& c' \frac{ (2k)! ^\frac{1}{2} }{(2k)^\frac{1}{4} (\sqrt{2t})^{2k}},\label{W2}
\end{IEEEeqnarray}
where we used again \eqref{AA13},
 and where $c$ and $c'$ denote some universal constants. The following lemma is needed.
\begin{lem}
\label{gevrey}
Let $s\in (0,1)$, and let $(a_k)_{k\ge 0} $ be a sequence  such that
\[
|a_k|\le C\frac{ k! ^s}{ R^k k^\alpha } \qquad \forall k\ge 0
\]
for some  constants $C>0,\ R>0$ and $\alpha \ge 0$.
Then the function
\[
f(x)=\sum_{k\ge 0} a_k \frac{x^k}{k!}
\]
is Gevrey of order $s$ on $[-L,L]$ for all $L>0$.
\end{lem}
\noindent
{\em Proof of Lemma \ref{gevrey}.} Note first that $f$ is well defined and analytic on $\mathbb R$, for
\[
|a_k|  \frac{L^k}{k!} \le C \sum_{k\ge 0} \frac{L^k}{R^k k^\alpha  k! ^{1-s}} < \infty .
\]
Taking the derivatives of the terms in the series, we infer that $f^{(m)}(x) = \sum_{k\ge 0 }a_{k+m} \frac{x^k}{k!}$, and hence
for any $x\in [-L,L]$
\begin{IEEEeqnarray}{rCl}
| f^{(m)} (x) | &\le& C\sum_{k\ge 0} \frac{L^k  (k+m) ! ^s}{ R^{k+m}  (k+m)^\alpha k!} \nonumber \\
&\le& C \frac{2^{sm} }{R^m m^\alpha } m!^s \sum_{k\ge 0} \frac{ (2^s L)^k  }{ R^k k! ^{1-s} }, \label{W1}
\end{IEEEeqnarray}
where we used the estimate
\begin{equation}
(k+m) ! \le 2^{k+m} k!\, m!.
\end{equation}
The proof of Lemma \ref{gevrey} is complete.\qed

It follows from \eqref{W2} and Lemma \ref{gevrey} that $e^{ i \frac{x^2}{4t} }$ is Gevrey of order $1/2$ in $x$. More precisely,
using \eqref{W2} and \eqref{W1}, we infer that for $x\in [-L,L]$ and $k\in \mathbb N$
\begin{equation}
|\partial _x ^k \big( e^{ i \frac{x^2}{4t} } \big) | \le C(t,L) \frac{k! ^\frac{1}{2} }{k^\frac{1}{4} (\sqrt{t} )^k } .
\end{equation}
Using \eqref{W3}, this yields for $x\in [-L,L]$ and $t_1\le t\le t_2$
\begin{IEEEeqnarray*}{rCl}
|\partial _x^k v(t,x) | &\le& C(t_1,t_2,L) \int_{-1}^1 \frac{k! ^\frac{1}{2} }{ k^\frac{1}{4} ( \sqrt{t})^k } |v_0(y)|dy\\
 &\le& C \frac{k! ^\frac{1}{2} }{ k^\frac{1}{4} (\sqrt{t})^k } ||\theta _0||_{L^1(0,1) }.
\end{IEEEeqnarray*}
Combined with \eqref{AA13} and \eqref{v1}, this gives for $(t,x)\in [t_1,t_2]\times [-L,L]$ and $(k,l)\in {\mathbb N} ^2$
\begin{IEEEeqnarray*}{rCl}
|\partial _x^k\partial_t^l  v(t,x) |
&=& |\partial _x^{k+2l}v(t,x)| \\
&\le&  C \frac{ (k+2l)! ^\frac{1}{2} }{ (k+2l)^\frac{1}{4} ( \sqrt{t} )^{k+2l}  }  ||\theta _0||_{L^1(0,1) }\\
&\le& C \frac{k!^\frac{1}{2} (2l)! ^\frac{1}{2}  2^{k+2l} }{ (k+2l) ^\frac{1}{4} (\sqrt{t}) ^{k+2l} } ||\theta _0||_{L^1(0,1) }  \\
&\le& C
\frac{k! ^\frac{1}{2} }{R_1^k} \frac{l!}{R_2^l} ||\theta _0||_{L^1(0,1) }
\end{IEEEeqnarray*}
where $C,R_1,R_2$ are some constants that depend on $t_1,t_2,$ and $L$.

Thus, $v$ is Gevrey of order 1 in $t$ and $1/2$ in $x$. In particular, $v(\tau , .)$ is an analytic function on $\mathbb R$. Being odd, it can be written as
\begin{equation}
v(\tau,x) =
\sum_{k\ge 0} y_k  (-i) ^k \frac{x^{2k+1 }}{ (2k +1 ) !}, \quad x\in \mathbb R \label{WWW1}
\end{equation}
with
\begin{IEEEeqnarray*}{rCl}
|y_k | &=& |\partial _x^{2k+1} v(\tau ,0) |\\
&\le& C (\tau )||\theta _0||_{ L^1(0,1) } \frac{ (2k+1) !^\frac{1}{2} }{(2k+1)^\frac{1}{4} (\sqrt{\tau}) ^{2k+1} } \\
&\le& C (\tau )||\theta _0||_{ L^1(0,1) }( \frac{2}{\tau } ) ^k  k ! ,
\end{IEEEeqnarray*}
where \eqref{AA13} was used again.
The proof of Lemma \ref{lem1} is complete. \qed

Since $v\in C({ \mathbb R}, L^2({\mathbb R} ))$, $\theta \in C([0,\tau ],L^2(0,1))$. \\

\noindent
{\sc Step 2. Construction of the control on $[\tau ,T]$}\\
We need the following
\begin{lem}
\label{lem2} Let $\tau \in (\frac{2}{3}T, T)$ and $1 < s < 2$ be numbers, and let $(y_k)_{k\ge 0}$ be as in \eqref{B12}.
Then there exists a function $y:[\tau , T]\to \mathbb R$ which is Gevrey of order $s$ on $[\tau , T ]$ and such that
\begin{IEEEeqnarray}{rCll}
y^{ (k) } (\tau )  &=& y_k \qquad &\forall k\ge 0,\label{Q1}\\
y^{ (k) } (T) &=& 0 \qquad &\forall k\ge 0, \label{Q2} \\
|y^{ (k) } (t)| &\le&  C ||\theta _0||_{L^1(0,1)}   \frac{ ( k ! )^s}{ R^k}, \ \ &\forall k\ge 0,\ \forall t\in [\tau , T] \quad \label{Q3}
\end{IEEEeqnarray}
for some constants $C=C(\tau ,T,s)>0$ and $R=R(\tau , T, s)>0$.
\end{lem}
\noindent
{\em Proof of Lemma \ref{lem2}.}  Let
\begin{equation}\label{eq:ybar}
\bar y(t) = \sum_{k \ge 0} y_k \frac{(t-\tau )^k}{k!} , \qquad t\in [\tau ,T].
\end{equation}
From $\tau >2T/3$, we infer that $2(T-\tau)/\tau <1$, and that for $t\in [\tau , T]$
\[
\sum_{k\ge 0} |y_k \frac{(t-\tau )^k}{k!} | \le C(\tau ) || \theta _0 ||_{ L^1(0,1) }  \sum_{k\ge 0} \left\vert \frac{ 2 (T-\tau) }{\tau } \right\vert ^k <\infty.
\]
Thus $\bar y$ is analytic on $[\tau , T]$, and hence Gevrey of order $s$ on $[\tau , T]$, with
\[
\bar y ^{ (k) } (\tau ) =  y_k, \qquad k\ge 0.
\]
Introduce the ``step function''
\[ \phi _s(t)
= \left\{
\begin{array}{ll}
1 & \text{ if } t\le 0,\\[2mm]
\displaystyle\frac{e^{ -(1-t )^{-\kappa } } }{ e^{ -(1-t)^{-\kappa} }   +  e^{- t^{-\kappa} } } &\text{ if } t\in (0,1),\\[3mm]
0 &\text{ if } t\ge 1,
\end{array}
\right.
\]
where $\kappa =(s-1)^{-1}$. Then $\phi _s$ is Gevrey of order $s$ on $[0,1]$ (in fact on $\mathbb R $) with
$\phi _s(0)=1$, $phi _s(1)=0$ and for all~$i\geq1$ $\phi _s^{ (i) } (0)=\phi _s^{ (i) } (1)=0$.

The desired function $y$ is given by
\[
y(t) = \phi _s \Bigl(\frac{t-\tau }{T-\tau }\Bigr) \bar y (t).
\]
Indeed, since products of Gevrey functions of order $s$ are Gevrey functions of order $s$, see e.g.
\cite{rudin,Yamanaka}, we infer that the function $y:[\tau ,T] \to\mathbb R$ is Gevrey of order $s$ on $[\tau , T]$. Furthermore,
\eqref{Q1} and \eqref{Q2} hold. Let us now prove \eqref{Q3}. Pick any $\rho >1$ such that  $2\rho (T-\tau) / \tau <1$. Then for any $z\in \mathbb C$ with
$|z|\le \rho (T-\tau)$, we have
\[
\sum_{k\ge 0} \Big|y_k \frac{z^k}{k!} \Big| \le M:= C(\tau ) || \theta _0 ||_{ L^1(0,1) }  \sum_{k\ge 0} \left\vert \frac{ 2 \rho (T-\tau) }{\tau } \right\vert ^k <\infty.
\]
It follows from Cauchy's formula that for $|z|\le T-\tau$ and $m\ge 0$,
\[
\left\vert \partial_z^m \sum_{k\ge 0} y_k \frac{z^k}{k!} \right \vert
\le M \frac{m !}{(\rho -1)^ m (T-\tau)^m }
\]
Thus, for $\tau \le t\le T$ and $m\ge 0$,
\begin{equation}
|{\bar y}^{(m)}(t)| \le M \frac{m !}{(\rho -1)^m(T-\tau ) ^m} \cdot  \label{QQ1}
\end{equation}
Then \eqref{Q3} follows from \eqref{QQ1} and
\cite[Theorem 19.7]{rudin}. The proof of Lemma \ref{lem2} is complete. \qed

Let, for $(t,x)\in [\tau ,T] \times [0,1] $,
\begin{IEEEeqnarray}{rCl}
\theta (t,x) &=& \sum_{k\ge 0} \frac{x^{2k +1 } }{ (2k +1) ! } (-i)^k y^{(k)} (t), \label{AA32}\\
u(t) &=& \sum _{k\ge 0} \frac{ (-i)^k y^{(k)} (t)}{(2k+1) !} \cdot  \label{AA31}
\end{IEEEeqnarray}

By Proposition \ref{prop1} and Lemma \ref{lem2}, the function $\theta$ is well-defined and Gevrey of order $s$ in $t$ and $s/2$ in $x$ on
$[\tau, T]\times [0,1]$. On the other hand, by Lemma \ref{lem1}, $\theta $ is Gevrey of order 1 in $t$ and $1/2$ in $x$ on $[ \varepsilon , \tau ] \times [0,1]$
for all $\varepsilon \in (0,\tau )$.
By \eqref{Q1}, the two series in \eqref{B110} and \eqref{AA32} take the same values at $t=\tau$,
so that $\theta\in C((0,T], C^k([0,1]) )$ for all $k\ge 0$ and $u\in C((0,T])$. Since $s>1$, to prove that $\theta$ is Gevrey of order $s$ in $t$
and $s/2$ in $x$ on $[\varepsilon , T]\times [0,1]$ for all $\varepsilon \in (0,T)$, it is sufficient
to notice  that for all $k\ge 0$ and $x\in [0,1]$
\begin{IEEEeqnarray*}{rCl}
(-i)^k\partial _t^k \theta (\tau ^-, x ) &=& \partial _x^{2k} \theta (\tau ^- ,x) \\
&=& \partial _x^{2k} \theta (\tau ^+ ,x)  \\
&=&  (-i)^k\partial _t^k \theta (\tau ^+, x).
\end{IEEEeqnarray*}
The proof of Theorem \ref{thm1} is complete. \qed

\end{pf}

\section{Controllability of the beam equation}
When $\theta$ is decomposed in terms of its real and imaginary parts, $\theta=\alpha+i\beta$, \eqref{B1} reads
\begin{IEEEeqnarray*}{rCl}
\alpha_t+\beta_{xx} &=& 0\\
\beta_t-\alpha_{xx} &=& 0;
\end{IEEEeqnarray*}
differentiating $\alpha$ w.r.t time and eliminating~$\beta$ yields
\begin{IEEEeqnarray*}{rCl}
\alpha_{tt}+\alpha_{4x} &=& 0,
\end{IEEEeqnarray*}
where $\alpha_{4x}:=\partial_x^4\alpha$. Hence we can adapt the results of the previous sections to derive the exact controllability of the Euler-Bernoulli beam equation
\begin{IEEEeqnarray}{rCl}
\eta _{tt} + \eta _{4x}  &=& 0\label{D1}\\
\bigl(\eta(t,0),\eta_{xx}(t,0)\bigr)&=&(0,0)\label{D2}\\
\bigl(\eta (t,1),\eta _{xx}(t,1)\bigr)&=&\bigl(u_1(t),u_2(t)\bigr)\label{D3} \\
\bigl(\eta(0,x),\eta _t(0,x)\bigr) &=& \bigl(\eta _0 (x),\eta _1(x)\bigr),\label{D4}
\end{IEEEeqnarray}
where $(t,x)\in (0,T)\times(0,1)$.

When $u_1=u_2=0$, \eqref{D1}-\eqref{D4} is a model for a simply supported
(or hinged) beam. It is well known, see e.g. \cite[Thm 6.16]{komornik-book}, that \eqref{D1}--\eqref{D4} is null controllable (or, equivalently, exactly controllable)
in $H^1_0(0,1)\times H^{-1}(0,1)$ in any time $T>0$ by using some controls $u_2\in L^2(0,T)$ and $u_1\equiv 0$.
The aim of this section is to derive in a straightforward way the exact controllability of \eqref{D1}-\eqref{D4} with {\em two} controls $u_1$ and $u_2$
by the flatness approach. There is no loss of generality in assuming that the final state is zero. The following result is a consequence of  Theorem \ref{thm1}.
\begin{thm}
\label{thm2}
Let $\eta _0\in H^2(0,1)\cap H^1_0(0,1)$, $\eta _1\in L^2(0,1)$ and $T>0$ be given.
Pick any $\tau \in (2T/3,T)$ and any $s\in (1,2)$. Then there exist $v_0\in H^2(\mathbb R)$ and  a function $y:[\tau , T]\to \mathbb R$
Gevrey of order $s$ on $[\tau , T]$ such that, setting
 \begin{equation}
u(t) = \left\{
\begin{array}{ll}
v(t,1)\quad &\text{if }\  0< t\le \tau ,\\
\sum _{k\ge 1} \frac{ (-i)^k y^{(k)} (t)}{(2k+1) !}  \quad &\text{if }\  \tau  < t \le T,
\end{array}
\right. \label{unew}
\end{equation}
(where $v$ denotes the solution of \eqref{v1}-\eqref{v2}) and $u_1(t) =\textrm{Re } u(t)$, $u_2(t)=\textrm{Im } u'(t)$,
the solution $\eta$ of \eqref{D1}-\eqref{D4} satisfies $\eta (T,.) =\eta_t (T,.)= 0$.  Furthermore, $u\in W^{1,4}(0,T)$ and it is Gevrey of order
$s$ in $t$ on $[ \varepsilon ,T]$ for all $\varepsilon \in (0,T)$,
$\eta \in C([0,T],L^2(0,1)) \cap C^\infty ((0,T]\times [0,1])$, and $\eta$ is Gevrey of order $s$ in $t$ and $s/2$ in $x$ on
$[\varepsilon ,T]\times [0,1]$ for all $\varepsilon \in (0,T)$.
\end{thm}
\begin{pf}
Let $\varphi=\eta _0$ and let $\psi$ solve $- \psi ''=\eta _1$, $\psi(0)=\psi (1)=0$. Then $\varphi$ and $\psi$ belong to $H^2(0,1)\cap H^1_0(0,1)$, and the same is true for
the (complex-valued) function $\theta _0(x):= \varphi(x)+i\psi (x)$. Pick a function $\zeta\in C^\infty (1,2)$ such that
\[
\zeta (x) =\left\{
\begin{array}{ll}
1 \quad &\textrm{if } x<5/4,\\
0 \quad &\textrm{if } x>7/4.
\end{array}
\right.
\]
Let $v_0:\mathbb R \to \mathbb C$ be defined as
\[
v_0 (x) =\left\{
\begin{array}{ll}
\theta _0(x) \quad &\textrm{if } 0<x<1,\\
-\theta _0(2-x) \zeta (x)  \quad &\textrm{if } 1<x<2,\\
0   \quad &\textrm{if } x>2,\\
-v_0(-x) \quad &\textrm{if } x<0.
\end{array}
\right.
\]
Then $v_0\in H^2(\mathbb R)$, $\textrm{supp }v_0\subset [ -7/4,7/4]$, and $v_0(-x)=-v_0(x)$ for all $x\in \mathbb R$. Clearly, $w_0=(v_0)_{xx} \in L^2(\mathbb R)$ satisfies also
$\textrm{supp }w_0\subset [ -7/4,7/4]$, and $w_0(-x)=-w_0(x)$ for a.e.  $x\in \mathbb R$. Let $v$ (resp. $w$) denote the solution of \eqref{v1}-\eqref{v2}
(resp. of  \eqref{v1}-\eqref{v2} with $w_0$ substituted to $v_0$). We know from \eqref{P1} that $v,w\in C^\infty ( {\mathbb R} _t \setminus \{0\} \times {\mathbb R} _x )$,
with  $w(t,x)=v_{xx}(t,x)$, and that
\[
\left( \int_{-\infty} ^\infty || v ||^4_{L^\infty ({\mathbb R})}  + || w ||^4_{L^\infty ({\mathbb R})}  dt \right)^{\frac{1}{4}} \le c
||v_0||_{H^2(\mathbb R  )}.
\]
Since $v_t=iv_{xx}=w$, we infer that
\[
v(.,1)\in W^{1,4}(0,\tau )\subset C([0,\tau ]).
\]
Let $\eta (t,x) =\textrm{Re } \theta (t,x)$, where $\theta _0$ is as above and $u$ is as in \eqref{unew}. Then the conclusion of Theorem \ref{thm1} (with the new function $v$)
is still valid, so that the regularity properties of $u$ and $\eta$ are established. Next,
\[
(\partial_t^2 + \partial _x^4) \theta =(-i\partial _t + \partial _x^2)(i\partial _t + \partial _x^2)\theta =0
\]
so that \eqref{D1} holds. \eqref{D2} (resp. \eqref{D3}) follows from \eqref{B1} and \eqref{B2} (resp. from \eqref{B1} and \eqref{B3}). \eqref{D4} is clear from
the construction of $\theta_0$, and we infer from Theorem \ref{thm1} and \eqref{B1} that $\eta(T,.)=\eta _t(T,.)=0$. \qed

\end{pf}

\section{Numerical experiments}
We illustrate the approach on a numerical example. The parameters are $\tau=0.35$, $T=0.5$ and~$s=1.9$; the (discontinuous) initial condition is
\begin{IEEEeqnarray*}{rCl}
    \textrm{Re~}\theta_0(x) &:=& \begin{cases}
    0& \text{if $x\in(0,0.5)$},\\
    1& \text{if $x\in(0.5,1)$}\end{cases}\\
     \textrm{Im~}\theta_0(x) &:=& \begin{cases}
    0& \text{if $x\in(0,0.2)\cup(0.7,1)$},\\
    1& \text{if $x\in(0.2,0.7)$}.\end{cases}   
\end{IEEEeqnarray*}

To compute the control $u(t)$ on~$[0,\tau]$ and the coefficients~$y_k$ we make use of the convolution formula~\eqref{C100}
\begin{IEEEeqnarray*}{rCl}
    v(t,x) &=& \int_{-1}^1 E(t,x-y)v_0(y)dy\\
    &=& \int_0^1 \underbrace{\bigl(E(t,x-y)-E(t,x+y)\bigr)}_{=:F(tx,y)}\theta_0(y)dy,
\end{IEEEeqnarray*}
where we have used~\eqref{v0}.
Then by~\eqref{eq:u}
\begin{IEEEeqnarray*}{rCl}
    u(t,x) = v(t,1) &=& \int_0^1 F(t,1,y)\theta_0(y)dy,
\end{IEEEeqnarray*}
and by~\eqref{B110}
\begin{IEEEeqnarray*}{rCl}
    y_k = \frac{\partial_x^{2k+1}v(\tau,0)}{(-i)^k} 
    &=&  \frac{1}{(-i)^k}\int_0^1 \partial_x^{2k+1}F(\tau,0,y)\theta_0(y)dy.
\end{IEEEeqnarray*}
All the integrals where numerically computed with the Matlab {\tt quadgk} function. The series~\eqref{eq:ybar} for $\bar y(t)$ was truncated at $k=15$, and so was the series~\eqref{AA31} for $u(t)$ on~$(\tau,T]$.

\begin{figure}[ht!]
\centering
\includegraphics[width=1.05\columnwidth]{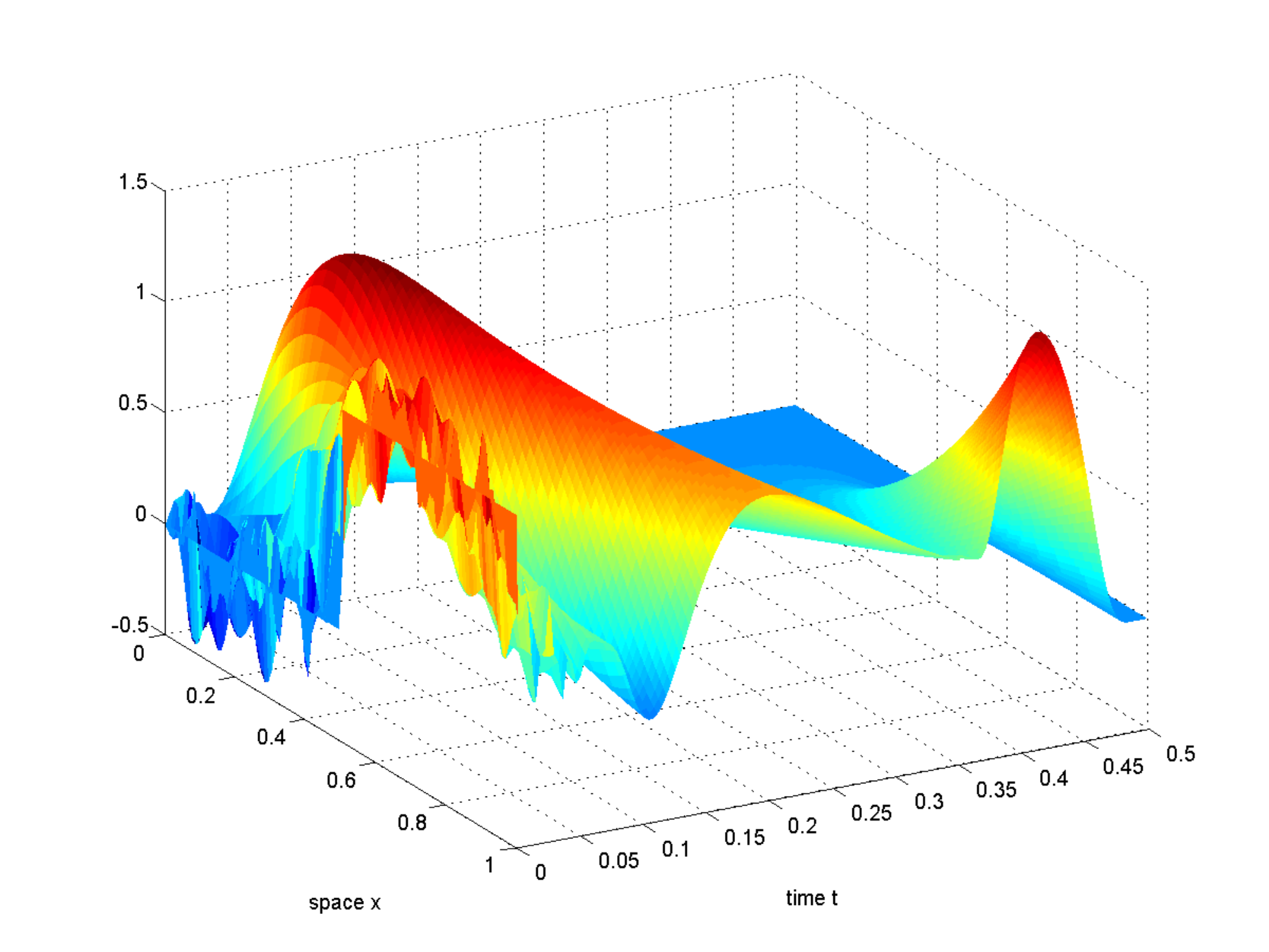}
\includegraphics[width=1.05\columnwidth]{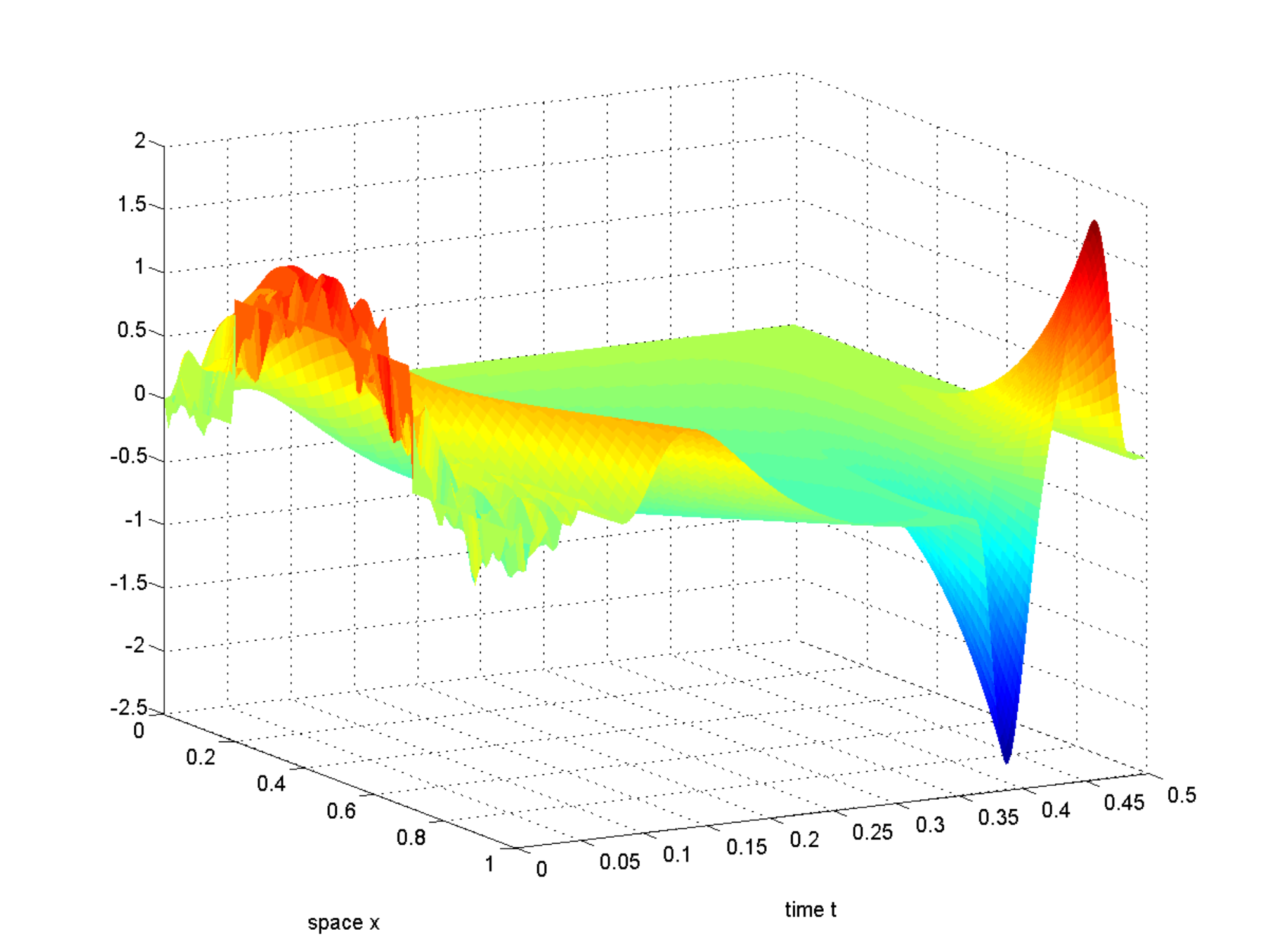}
\caption{$\textrm{Re~}\theta$ (top) and $\textrm{Im~}\theta$ (bottom).}
\label{fig:surtheta}
\end{figure}

Fig.~\ref{fig:surtheta} displays the evolution of~$\theta$ on~$[0,T]$. The regularizing effect of the control on~$(0,\tau)$ is clearly visible.

\bibliography{schro1D-BIB}

\end{document}